\documentclass[a4paper,10pt]{amsproc}
\usepackage{amssymb}
\usepackage{graphicx}
\usepackage{amscd}
\usepackage{amsmath}

 \theoremstyle{plain}

\newtheorem{lemma}{Lemma}

\newtheorem{proposition}{Proposition}
\newtheorem{remark}{Remark}

\newtheorem{theorem}{Theorem}

\numberwithin{equation}{section}

\begin{document}

\title[Inverse problems for finite Jacobi matrices and
Krein--Stieltjes strings.]{Inverse problems for finite Jacobi
matrices and Krein--Stieltjes strings.}

\author{Alexander Mikhaylov} 
\address{St. Petersburg   Department   of   V.A. Steklov    Institute   of   Mathematics
of   the   Russian   Academy   of   Sciences, 7, Fontanka, 191023
St. Petersburg, Russia and Saint Petersburg State University,
St.Petersburg State University, 7/9 Universitetskaya nab., St.
Petersburg, 199034 Russia.} \email{mikhaylov@pdmi.ras.ru}

\author{Victor Mikhaylov} 
\address{St.Petersburg   Department   of   V.A.Steklov    Institute   of   Mathematics
of   the   Russian   Academy   of   Sciences, 7, Fontanka, 191023
St. Petersburg, Russia} \email{vsmikhaylov@pdmi.ras.ru}


\maketitle

\begin{abstract}
We consider dynamic inverse problems for a dynamical system
associated with a finite Jacobi matrix and for a system describing
propagation of waves in a finite Krein-Stieltjes string. We offer
three methods of recovering unknown parameters: entries of a
Jacobi matrix in the first problem and point masses and distances
between them in the second, from dynamic Dirichlet-to-Neumann
operators. We also answer a question on a characterization of
dynamic inverse data for these two problems.

\end{abstract}

\section{Introduction}

Given a sequence of positive numbers $\{a_0,$
$a_1,\ldots,a_{N-1}\}$ (in what follows we assume $a_0=1$) and
real numbers $\{b_1, b_2,\ldots,b_N \}$, we denote by $A$ the
finite Jacobi matrix given by
\begin{equation}
\label{Jac_matr}
A=\begin{pmatrix} b_1 & a_1 & 0 & 0 & 0 &\ldots \\
a_1 & b_2 & a_2 & 0 & 0 &\ldots \\
\ldots &\ldots  &\ldots &\ldots & \ldots &\ldots \\
\ldots &\ldots  &\ldots &a_{N-2} & b_{N-1} &a_{N-1}\\
\ldots &\ldots &\ldots &\ldots & a_{N-1} &b_N
\end{pmatrix}.
\end{equation}
Let $u=(u_1,\ldots,u_N)\in \mathbb{R}^N$ and $T>0$ be fixed. With
the matrix $A$ we associate the dynamical system:
\begin{equation}
\label{DnSst}
\begin{cases}
u_{tt}(t)-Au(t)=F(t),\quad t>0,\\
u(0)=u_t(0)=0,
\end{cases}
\end{equation}
where the vector function $F(t)=(f(t),0,\ldots,0)$, $f\in
L_2(0,T)$ is interpreted as a \emph{boundary control}. The
solution of (\ref{DnSst}) is denoted by $u^f$. With the system
(\ref{DnSst}) we associate the \emph{response operator} acting by
the rule
\begin{equation}
\label{RespJM} \left(R^Tf\right)(t)=u^f_1(t),\quad 0<t<T.
\end{equation}
In the first part of the paper we will be dealing with the dynamic
inverse problem (IP) of the reconstruction of the matrix $A$ from
the response operator $R^T$, in the second part we study the
related IP for a finite Krein--Stieltjes string.

The forward and inverse spectral problems for finite, infinite and
periodic Jacobi operators are subjects of numerous papers, to
mention just \cite{T} and references therein. At the same time,
corresponding dynamic problems were not studied in the literature
except for the author's papers \cite{MM,MM1,MM2} on inverse
problems for dynamical system associated with Jacobi matrix in
which the time is assumed to be discrete; in \cite{MMS,MM4} this
model is applied to the studying of a Weyl function and moment
problems; and \cite{MM3} where the model with continuous time for
a Krein-Stieltjes string was considered.

In \cite{MM1,MM11} the authors established the relationships
between de Branges method \cite{DBr,DMcK} for IP for canonical
systems and the Boundary Control (BC) method for dynamic IP
\cite{B07,B17}. The common strategy in studying (spectral) IP for
general canonical systems is to approximate the Hamiltonian by the
sequence of finite-rank matrices of special type, see
\cite{DBr,DMcK,A}. That is why we think that the present paper,
devoted to the most simple case of the dynamic IP for the
finite-dimensional dynamical system, can be considered as a first
step toward the solution of a general IP for dynamical systems
associated with semi-infinite Jacobi matrices and for general
Krein--Stieltjes strings.

In the second section we focus on the IP for dynamical system
(\ref{DnSst}). We derive the representation for the solution to
(\ref{DnSst}), introduce operators of the BC method
\cite{B07,B17}, propose three methods of solving the IP and give a
characterization of dynamic inverse data in terms of properties of
so-called connecting operator. The last section is devoted to the
special case of dynamical system (\ref{DnSst}): Krein--Stieltjes
string. We briefly outline the set up of the problem following the
brief note \cite{MM3}, since some of results formulated in
\cite{MM3} were incorrect, we focus on the derivation of Krein
equations and results on the characterization of inverse data.

\section{IP for dynamical system associated with finite Jacobi matrix.}

\subsection{Forward problem}
The following Cauchy problem for the difference equation
\begin{equation*}
\begin{cases}
a_1\phi_2+b_1\phi_1=\lambda\phi_1,\\
a_n\phi_{n+1}+a_{n-1}\phi_{n-1}+b_n\phi_n=\lambda \phi_n,\quad
n=2,\ldots,N,\\
\phi_1=1,
\end{cases}
\end{equation*}
determines the set of polynomials
$\{1,\phi_2(\lambda),\ldots,\phi_N(\lambda),\phi_{N+1}(\lambda)\}$.
Denote by $\lambda_1,\ldots,\lambda_{N}$ roots of the equation
$\phi_{N+1}(\lambda)=0$. It is known \cite{A,Ah} that they are
real and distinct. By $(\cdot,\cdot)$ we denote the scalar product
in $\mathbb{R}^N$ and introduce vectors and coefficients by the
rules:
\begin{equation*}
\varphi(\lambda)=\begin{pmatrix} \phi_1(\lambda)\\
\cdot\\ \cdot\\ \phi_{N}(\lambda)\end{pmatrix},\quad\varphi_k=\begin{pmatrix} \phi_1(\lambda_k)\\
\cdot\\ \cdot\\ \phi_{N}(\lambda_k)\end{pmatrix},\quad
\rho_k=\left(\varphi_k,\varphi_k\right),\quad k=1,\ldots,N.
\end{equation*}
Thus $\varphi_k$ are non-normalized eigenvectors of $A$,
corresponding to eigenvalues $\lambda_k$:
\begin{equation*}
A\varphi_k=\lambda_k\varphi_k,\quad \quad k=1,\ldots,N.
\end{equation*}
We call by \emph{spectral data} and spectral function $\rho$ the
following objects:
\begin{equation}
\label{measure_JM} \left\{\lambda_i,\rho_i\right\}_{i=1}^{N},\quad
\rho(\lambda)=\sum_{\{k:\lambda_k<\lambda\}}\frac{1}{\rho_k}.
\end{equation}
Note the property:
\begin{equation}
\label{rho_pr} \sum_{k=1}^N\frac{1}{\rho_k}=1.
\end{equation}

The standard application of a Fourier method (see
\cite{AI,AM,MM1}) yields
\begin{lemma}
The solution to (\ref{DnSst}) admits the spectral representation
\begin{equation}
\label{Sol_spec_repr_JM}
u^f(t)=\sum_{k=1}^{N}h_k(t)\varphi_k,\quad
u^f(t)=\int_{-\infty}^\infty\int_0^t
S(t-\tau,\lambda)f(\tau)\,d\tau\varphi(\lambda)\,d\rho(\lambda).
\end{equation}
where
\begin{eqnarray}
h_k(t)=\frac{1}{\rho_k}\int_0^t
f(\tau)S_k(t-\tau)\,d\tau,\notag\\
S(t,\lambda)=\begin{cases}
\frac{\sin{\sqrt{\lambda}t}}{\sqrt{\lambda}},\quad
\lambda>0,\\
\frac{\operatorname{sh}{\sqrt{|\lambda|}t}}{\sqrt{|\lambda|}},\quad
\lambda<0,\\
t,\quad \lambda=0,
\end{cases},\quad S_k(t)=S(t,\lambda_k).\label{S_k}
\end{eqnarray}
\end{lemma}

\subsection{Operators of the BC method.}

We introduce the \emph{outer space} of the system (\ref{DnSst}),
the space of controls: $\mathcal{F}^T:=L_2(0,T)$ with the scalar
product $f,g\in \mathcal{F}^T,$
$\left(f,g\right)_{\mathcal{F}^T}=\int_0^Tf(t)g(t)\,dt$. The
\emph{response operator} $R^T: \mathcal{F}^T\mapsto \mathcal{F}^T$
is introduced by formula (\ref{RespJM}). Making use of
(\ref{Sol_spec_repr_JM}) implies the representation for $R^T$:
\begin{equation}
\label{Resp_def}
\left(R^Tf\right)(t)=u_1^f(t)=\sum_{k=1}^{N}h_k(t)=\int_0^t
r(t-s)f(s)\,ds,
\end{equation}
where
\begin{equation}
\label{resp_func_JM}
r(t)=\sum_{k=1}^{N}\frac{1}{\rho_k}S_k(t),
\end{equation}
is called a \emph{response function}. Note that the operator $R^T$
is a natural analog of a dynamic Dirichlet-to-Neumann operator
\cite{B07,B17} in continuous, and \cite{MM,MM2,MMS} in discrete
cases.

The \emph{inner space}, i.e. the space of states of (\ref{DnSst})
is denoted by $\mathcal{H}^N:=\mathbb{R}^{N}$, indeed for any
$T=N>0$ we have that $u^f(T)\in \mathcal{H}^N$. The metric in
$\mathcal{H}^N$ is given by
\begin{equation*}
(a,b)_{\mathcal{H}^N}=\left(a,b\right)_{\mathbb{R}^{N}},\quad
a,b\in\mathcal{H}^N.
\end{equation*}
The \emph{control operator} $W^T: \mathcal{F}^T\mapsto
\mathcal{H}^N$ is introduced by the rule:
\begin{equation*}
W^Tf=u^f(T).
\end{equation*}
Due to (\ref{Sol_spec_repr_JM}) we have that
$W^Tf=\sum_{k=1}^{N}h_k(T)\varphi_k$. The boundary controllability
properties of a dynamical system plays a crucial role in a
procedure of solving the dynamic inverse problems, see
\cite{B17,MM2,BM_1,AM,B08,MM1}. We set up the following
\emph{control problem}: for a fixed state $a\in \mathcal{H}^N$,
$a=\sum_{k=1}^{N}a_k\varphi_k$ we look for a control $f\in
\mathcal{F}^T$ driving the system (\ref{DnSst}) to a prescribed
state:
\begin{equation}
\label{control_JM} W^Tf=a.
\end{equation}
Using (\ref{Sol_spec_repr_JM}) we see that the equality
(\ref{control_JM}) is equavalent to the following moment problem:
to find $f\in \mathcal{F}^T$ such that
\begin{equation*}
a_k=h_k(T)=\frac{1}{\rho_k}\int_0^T f(\tau)S_k(T-\tau)d\tau,\quad
k=1,\ldots,N.
\end{equation*}
Clearly, a solution to this moment problem exists, but it is not
unique. We introduce the subspace
\begin{equation*}
\mathcal{F}^T_1={\operatorname{span}\left\{S_k(T-t)\right\}_{k=1}^N}.
\end{equation*}
The following lemma (see \cite{AI}) answers the question on the
controllability of (\ref{DnSst}):
\begin{lemma}
\label{LemmaCont_JM} The operator ${W}^T$ maps $\mathcal{F}^T_1$
onto $\mathcal{H}^N$ isomorphically.
\end{lemma}
\begin{proof}
Let $\left\{\widetilde S_k(T-t)\right\}_{k=1}^N,$ where
$\widetilde S_k(T-t)\in \mathcal{F}^T_1$ be a bi-orthogonal basis
in $\mathcal{F}^T_1$, i.e. $\int_0^T S_k(T-t)\widetilde
S_l(T-t)\,dt=\delta_{kl}.$ Then we immediately have that
\begin{equation*}
f(t)=\sum_{k=1}^N a_k\rho_k\widetilde S_k(T-t).
\end{equation*}
\end{proof}

The \emph{connecting operator} $C^T:\mathcal{F}^T\mapsto
\mathcal{F}^T$ is defined by the rule
$C^T:=\left(W^T\right)^*W^T$, so by the definition for $f,g\in
\mathcal{F}^T$ one has
\begin{equation}
\label{C_T_JM}
\left(C^Tf,g\right)_{\mathcal{F}^T}=\left(u^f(T),u^g(T)\right)_{\mathcal{H}^N}=\left(W^Tf,W^Tg\right)_{\mathcal{H}^N}.
\end{equation}
It is crucial in the BC method that $C^T$ can be expressed in
terms of inverse data:
\begin{theorem}
The connecting operator admits the representation in terms of
dynamic inverse data:
\begin{equation}
\label{CT_dyn_repr_JM}
\left(C^Tf\right)(t)=\frac{1}{2}\int_0^T\int_{|t-s|}^{2T-s-t}r(\tau)\,d\tau
f(s)\,ds,
\end{equation}
and in terms of spectral inverse data:
\begin{equation}
\label{CT_sp_repr_JM} \left(C^Tf\right)(t)=\int_0^T
\sum_{k=1}^{N}\frac{1}{\rho_k}S_k(T-t)S_k(T-s)f(s)\,ds.
\end{equation}
\end{theorem}
\begin{proof}
Taking arbitrary $f,g\in \mathcal{F}^T$ and introducing the
function $\psi(t,s):=\left(u^f(t),u^g(s)\right)_{\mathcal{H}^N}$,
one can show that $\psi$ satisfies the equation
\begin{equation*}
\psi_{tt}-\psi_{ss}=\frac{1}{l_1}\left(f(t)(Rg)(t)-g(s)(Rf)(t)\right).
\end{equation*}
Solving it by the d'Alembert method and noting that
$\psi(T,T):=\left(C^Tf,g\right)_{\mathcal{F}^T}$, yields
(\ref{CT_dyn_repr_JM}). The formula (\ref{CT_sp_repr_JM}) follows
from (\ref{C_T_JM}) and (\ref{Sol_spec_repr_JM}).
\end{proof}

\begin{remark}
The formula (\ref{CT_sp_repr_JM}) implies that
$\mathcal{F}^T_1=C^T\mathcal{F}^T$, that is, $\mathcal{F}^T_1$ is
completely determined by inverse data. The important properties of
functions from $\mathcal{F}^T_1$ is that
\begin{equation}
\label{FT_proper_JM} f(T)=0,\quad \text{for all}\,\, f\in
\mathcal{F}^T_1.
\end{equation}
\end{remark}

\subsection{Inverse problem.}

Using spectral theorem we can write down the solution to
(\ref{DnSst}) in the form:
\begin{equation}
\label{U_repr}
u^f(t)=\int_0^tA^{-\frac{1}{2}}\sin\left[A^{\frac{1}{2}}(t-s)\right]F(s)\,ds.
\end{equation}
Here and everywhere below for $A^{\frac{1}{2}}$ we choose the
positive branch of square root.  The definition (\ref{Resp_def})
and (\ref{U_repr}) imply that the response function is given by
the following matrix element:
\begin{equation*}
r(t)=u^\delta_1(t)=\left(A^{-\frac{1}{2}}\sin\left[
{A}^{\frac{1}{2}}t\right]\right)_{11},
\end{equation*}
where $\delta(t)$ is a Dirac delta-function. Expanding the above
equality in Tailor series at zero we get that
\begin{equation}
\label{Tailor1}
r(t)=t-\frac{1}{3!}\left(A\right)_{11}t^3+\frac{1}{5!}\left(A^2\right)_{11}t^5-\frac{1}{7!}\left(A^3\right)_{11}t^7+\ldots,\quad
t\to 0.
\end{equation}
On the other hand the spectral representation
(\ref{Sol_spec_repr_JM}) implies the following relation:
\begin{equation}
\label{Tailor2} r(t)=\int_{\mathbb{R}}
S(t,\lambda)\,d\rho(\lambda).
\end{equation}
Equating the coefficients in Taylor expansion of $r$ at zero in
(\ref{Tailor1}) and (\ref{Tailor2}) yields the following relations
between derivatives of response function at zero, coefficients of
Jacobi matrix $A$ and moments of spectral measure of $A$ (we write
down first several terms):
\begin{eqnarray*}
r'(0)=1=\int_{\mathbb{R}}1\,d\rho(\lambda)=s_0,\\
r^{(3)}(0)=A_{1\,,1}=b_1=-\int_{\mathbb{R}}\lambda\,d\rho(\lambda)=s_1,\\
r^{(5)}(0)=\left(A^2\right)_{1\,1}=b_1^2+a_1^2=\int_{\mathbb{R}}\lambda^2\,d\rho(\lambda)=s_2,\\
r^{(7)}(0)=\left(A^3\right)_{1\,1}=b_1^3+2b_1a_1^2+a_1b_2=-\int_{\mathbb{R}}\lambda^3\,d\rho(\lambda)=s_3.
\end{eqnarray*}
Note that these trace-type formulas allow one to recover
coefficients $a_k,b_k$, $k=1,2,\ldots$ recursively from $r(t)$
given on an interval $(0,\varepsilon)$ for any $\varepsilon>0$.

We introduce the subspace $\mathcal{F}_{0}^T:=\{f\in
C^\infty(0,T):\ f(0)=f(T)=f'(0)=f'(T)=0 \}$. The representation
(\ref{Sol_spec_repr_JM}) implies that for $f\in \mathcal{F}^T_0$
one has that $u^{f_{tt}}=u_{tt}^f$. We use this fact taking
$f,g\in \mathcal{F}_{0}^T$ and evaluating the following quadratic
form:
\begin{eqnarray}
\left(C^Tf_{tt},g\right)_{\mathcal{F}^T}=\left(u^f_{tt}(T),u^g(T)\right)_{\mathbb{R}^{N}}=
\left(u^f_{tt}(T),u^g(T)\right)_\mathcal{H^N} \notag\\
=\left(Au^f(T)+F(T),u^g(T)\right)_\mathcal{H^N}
=\left(Au^f(T),u^g(T)\right)_\mathcal{H^N}.\label{L_quadr_JM}
\end{eqnarray}
Then it is possible to perform the spectral analysis of $A$ using
the classical variational approach, controllability (Lemma
\ref{LemmaCont}) of the system (\ref{DnSst}) and representation
(\ref{L_quadr_JM}), see also \cite{B07,B01JII}. Then the spectral
data can be recovered by the following procedure:
\begin{itemize}
\item[1)] The first eigenvalue is given by
\begin{equation}
\label{M1}\lambda_1=\inf_{f\in \mathcal{F}^T_0,\,(C^T
f,f)_{\mathcal{F}^T}=1}\left(C^T f_{tt},f\right)_{\mathcal{F}^T}.
\end{equation}

\item[2)] Let $\left\{f^1_k\right\}_{k=1}^\infty$, be any
minimizing sequence of (\ref{M1}), then
\begin{equation*}
\rho_1=\frac{1}{\left(\varkappa_1\right)^2},\quad \text{where
$\varkappa_1=\lim_{k\to\infty}\left(R^Tf^1_k\right)(T)$}.
\end{equation*}

\item[3)] The second eigenvalue is given by
\begin{equation}
\label{M2}\lambda_2=\inf_{\substack {f\in
\mathcal{F}^T_0,(C^T f,f)_{\mathcal{F}^T}=1\\
\lim_{k\to\infty}(C^T f,f_k^1)_{\mathcal{F}^T}=0}}\left(C^T
f_{tt},f\right)_{\mathcal{F}^T},
\end{equation}
where  $\left\{f^1_k\right\}_{k=1}^\infty$ is any minimizing
sequence of (\ref{M1}).

\item[4)] Let $\left\{f^2_k\right\}_{k=1}^\infty$, be the
minimizer of (\ref{M2}), then
\begin{equation*}
\rho_2=\frac{1}{\left(\varkappa_2\right)^2},\quad \text{where
$\varkappa_2=\lim_{k\to\infty}\left(R^Tf^2_k\right)(T)$}.
\end{equation*}
\end{itemize}
We note that the subspace $\mathcal{F}^T_0$ consists of functions
with additional conditions on the boundary, but these conditions
does not necessarily hold for minimizer. Continuing the described
procedure, we recover the set $\{\lambda_k,\rho_k\}_{k=1}^{N}$ and
construct the measure $d\rho(\lambda)$ by (\ref{measure_JM}).
After that the matrix $A$ can be recovered by known spectral
methods, see for example \cite{MM2}.

\subsection{Krein equations.}
\label{Kr_eq}

By $f_k^T\in \mathcal{F}^T_1$ we denote controls, driving the
system (\ref{DnSst}) to prescribed \emph{special states}
\begin{equation*}
d_k\in \mathcal{H}^N, \,d_k=\left(0,\ldots,1,\ldots,0\right),\quad
k=1,\ldots,N.
\end{equation*}
In other words, $f^T_k$ are the solutions to the equations
$W^Tf_k^T=d_k$, $k=1\,\ldots,N$. By Lemma \ref{LemmaCont_JM} we
know that such controls exist and are unique. It is important that
they can be found as the solutions to the Krein equations.

\begin{theorem} The control $f_1^T$ can be found as the solution to
the following equation:
\begin{equation}
\label{Cont_f1_JM} \left(C^Tf_1^T\right)(t)=r(T-t),\quad 0<t<T.
\end{equation}
The controls $f^T_k$ satisfy the system:
\begin{equation}
\label{Cont_system_JM}
\begin{cases}
-\left(C^Tf_1^T\right)''=b_1C^Tf^T_1+a_1C^Tf^T_{2},\\
-\left(C^Tf_k^T\right)''=a_{k-1}C^Tf^T_{k-1}+b_kC^Tf^T_k+a_kC^Tf^T_{k+1},\quad
k=2,\ldots,N-1,\\
-\left(C^Tf_N^T\right)''=a_{N-1}C^Tf^T_{N-1}+b_NC^Tf^T_N,
\end{cases}
\end{equation}
\end{theorem}
\begin{proof}
Let us check (\ref{Cont_f1_JM}): for $g\in \mathcal{F}^T$ we have:
\begin{equation*}
\left(C^Tf_1^T,g\right)_{\mathcal{F}^T}=\left(u^{f_1^T}(T),u^g(T)\right)_{\mathcal{H}^N}=u_1^g(T)=\int_0^Tr(T-\tau)g(\tau)\,d\tau,
\end{equation*}
which, due to the arbitrariness of $g$, yields (\ref{Cont_f1_JM}).

Taking some $k=1,\ldots,N$, we rewrite (\ref{Sol_spec_repr_JM})
with $f_k^T$, then for $n-$th component of $u^{f_k^T}$ we have the
representation
\begin{equation*}
u^{f_k^T}_n(T)=\sum_{i=1}^{N}\frac{1}{\rho_i}\int_0^T
f^T_k(\tau)S_i(T-\tau)\,d\tau\varphi_i^n=\begin{cases}
1,\quad n=k,\\
0,\quad n\not=k.
\end{cases}
\end{equation*}
multiplying the above equality by $\varphi_j^n$ and summing up
from $n=1$ to $N$, using the orthogonality property of $\varphi_k$
\cite{Ah,A} we get that
\begin{equation}
\label{f1} \int_0^T
f_k^T(\tau)S_j(T-\tau)\,d\tau=\varphi_j^k,\quad j=1,\ldots,N.
\end{equation}
Using (\ref{CT_sp_repr_JM}) and (\ref{f1}) we obtain that
\begin{equation}
\label{f2}
\left(C^Tf^T_k\right)(t)=\sum_{n=1}^N\frac{1}{\rho_n}\varphi_n^k
S_n(T-t).
\end{equation}
Differentiating twice (\ref{f2}) and counting that
\begin{align*}
&\lambda_n\varphi_n^1=b_1\varphi_n^1+a_1\varphi_n^2,\\
&\lambda_n\varphi_n^k=b_k\varphi_n^k+a_{k-1}\varphi_n^{k-1}+a_k\varphi_n^{k+1},\quad
k=2,\ldots,N-1,\\
&\lambda_n\varphi_n^N=a_{N-1}\varphi_n^{N-1}+b_N\varphi_n^N,
\end{align*}
we get:
\begin{align*}
&\left(C^Tf^T_1\right)''(t)=\sum_{n=1}^N\frac{-\lambda_n}{\rho_n}\varphi_n^1S_n(T-t)
=-\sum_{n=1}^N\frac{b_1\varphi_n^1+a_1\varphi_n^{2}}{\rho_n}S_n(T-t)\\
&=-\left(b_1\left(C^Tf^T_1\right)(t)+a_1\left(C^Tf^T_{2}\right)(t)\right),\\
&\left(C^Tf^T_k\right)''(t)=\sum_{n=1}^N\frac{-\lambda_n}{\rho_n}\varphi_n^kS_n(T-t)
=-\sum_{n=1}^N\frac{b_k\varphi_n^k+a_{k-1}\varphi_n^{k-1}+a_k\varphi_n^{k+1}}{\rho_n}S_n(T-t)\\
&=-\left(b_k\left(C^Tf^T_k\right)(t)+a_{k-1}\left(C^Tf^T_{k-1}\right)(t)+a_k\left(C^Tf^T_{k+1}\right)(t)\right),\quad
k=2,\ldots,N-1,\\
&\left(C^Tf^T_N\right)''(t)=\sum_{n=1}^N\frac{-\lambda_n}{\rho_n}\varphi_n^NS_n(T-t)
=-\sum_{n=1}^N\frac{a_{N-1}\varphi_n^{N-1}+b_N\varphi_n^{N}}{\rho_n}S_n(T-t)\\
&=-\left(a_{N-1}\left(C^Tf^T_{N-1}\right)(t)+b_N\left(C^Tf^T_{N}\right)(t)\right),
\end{align*}
which yields (\ref{Cont_system_JM}).


\end{proof}

To recover unknown coefficient of $A$ we propose the following
procedure. First we note that by the definition of controls
$f^T_k$:
\begin{equation}
\label{F_orthog}
\left(C^Tf^T_i,f^T_j\right)_{\mathcal{F}^T}=\delta_{ij}.
\end{equation}
Thus, multiplying the first line in (\ref{Cont_system_JM}) by
$f^T_1$ in $\mathcal{F}^T$, we obtain that
\begin{equation*}
b_1=-\left(\left(C^Tf^T_1\right)'',f^T_1\right)_{\mathcal{F}^T}.
\end{equation*}
Then from the first line in (\ref{Cont_system_JM}) we have that
\begin{equation*}
a_1C^Tf^T_{2}=-\left(C^Tf_1^T\right)''-b_1C^Tf^T_1.
\end{equation*}
Inverting $C^T$ in $\mathcal{F}^T_1$, we can find the quantity
$a_1f^T_{2}$. Then evaluating the quadratic form
\begin{equation*}
\left(C^Ta_1f^T_2,a_1f^T_2\right)_{\mathcal{F}^T}=a_1^2\left(C^Tf^T_2,f^T_2\right)_{\mathcal{F}^T}=a_1^2,
\end{equation*}
we recover $a_1$ and $f^T_2$. Continuing this procedure we obtain
that
\begin{eqnarray*}
b_n=-\left(\left(C^Tf^T_n\right)'',f^T_n\right)_{\mathcal{F}^T},\\
a_n=-\left(\left(C^Tf^T_n\right)'',f^T_{n+1}\right)_{\mathcal{F}^T}=-\left(\left(C^Tf^T_{n+1}\right)'',f^T_n\right)_{\mathcal{F}^T}.
\end{eqnarray*}


Before we turn to the characterization of inverse data, we
formulate three propositions:
\begin{proposition}
\label{Prop_contr} Components of eigenvector $\varphi_n$ admit the
following representations:
\begin{equation}
\varphi^k_n=\int_0^Tf^T_k(\tau)S_n(T-\tau)\,d\tau,\quad
k=1,\ldots,N.
\end{equation}
\end{proposition}
\begin{proof}
Writing down the representation
\begin{equation*}
u^{f^T_k}(T)=\sum_{l=1}^{N}\frac{1}{\rho_l}\int_0^T
f^T_k(\tau)S_l(t-\tau)\,d\tau \varphi_l=\left(0,\ldots,1\ldots,
0\right)^T
\end{equation*}
with $1$ being at $k-$th place we multiply the above equality by
$\varphi_n$ to get:
\begin{equation*}
\varphi_n^k=\int_0^T f^T_k(\tau)S_n(t-\tau)\,d\tau,
\end{equation*}
which proves the statement.
\end{proof}

\begin{proposition}
\label{Prop2} If the response function is given by
(\ref{resp_func_JM}), then for the connecting operator defined by
(\ref{CT_dyn_repr_JM}) the following relation holds:
\begin{eqnarray*}
\left(\left(C^Tf\right)'',g\right)=\left(f,\left(C^Tg\right)''\right)\\
=\frac{1}{2}\int_0^T\int_0^T\left(r'(2T-t-s)-r'(|t-s|)\right)f(s)g(t)\,ds\,dt,
\end{eqnarray*}
for $f,g\in C^1(0,T)$.
\end{proposition}
\begin{proof}
We rewrite (\ref{CT_dyn_repr_JM}) in the form:
\begin{equation*}
2\left(C^Tf\right)(t)=\int_0^t\int_0^{2T-s-t}r(\tau)\,d\tau
f(s)\,ds+\int_t^T\int_{s-t}^{2T-s-t}r(\tau)\,d\tau f(s)\,ds
\end{equation*}
Then differentiating we obtain
\begin{eqnarray*}
2\left(C^Tf\right)'(t)=\int_0^{2T-2t}r(\tau)\,d\tau+\int_0^t
f(s)[-r(2T-s-t)-r(t-s)]\,ds\\
-\int_0^{2T-2t}r(\tau)\,d\tau+ \int_t^T
f(s)[-r(2T-s-t)+r(s-t)]\,ds\\ =\int_0^T
f(s)[r(s-t)-r(2T-s-t)]\,ds,
\end{eqnarray*}
where from (\ref{resp_func_JM}) we used the fact that $r$ is odd.
Taking the second derivative we get
\begin{equation*}
2\left(C^Tf\right)''(t)=\int_0^T f(s)[r'(2T-s-t)-r'(|t-s|)]\,ds,
\end{equation*}
which immediately yields the statement.
\end{proof}

\begin{proposition}
\label{Prop3} If the response function is given by
(\ref{resp_func_JM}), norming coefficients satisfy (\ref{rho_pr})
and the connecting operator is defined by (\ref{CT_dyn_repr_JM}),
(\ref{CT_sp_repr_JM}), then if $f\in \mathcal{F}^T$ satisfies the
equation
\begin{equation}
\label{Eq1} \left(C^Tf\right)(t)=r(T-t),
\end{equation}
then for such $f$ the equality $\left(C^Tf,f\right)=1$ holds.
\end{proposition}
\begin{proof}
We rewrite (\ref{CT_sp_repr_JM}) in the form
\begin{equation*}
\left(C^Tf\right)(t)=\sum_{k=1}^{N}f_k(T)S_k(T-t),\quad
f_k(T)=\frac{1}{\rho_k}\int_0^TS_k(T-s)f(s)\,ds.
\end{equation*}
Then the representation of response function (\ref{resp_func_JM})
and equation (\ref{Eq1}) implies that for such $f$ the following
relation holds:
\begin{equation*}
f_k(T)=\frac{1}{\rho_k},\quad k=1,\ldots,N.
\end{equation*}
Then we can evaluate
\begin{eqnarray*}
\left(C^Tf,f\right)=\left(r(T-\cdot),f\right)=\int_0^T\sum_{k=1}^{N}\frac{1}{\rho_k}f(t)S_k(T-t)\,dt\\
=\sum_{k=1}^{N}f_k(T)=\sum_{k=1}^{N}\frac{1}{\rho_k}=1.
\end{eqnarray*}

\end{proof}

\begin{theorem}\label{thm3}
The function $r$ is a response function for the dynamical system
(\ref{DnSst}) if and only if it has a form (\ref{resp_func_JM})
where coefficients satisfy (\ref{rho_pr}), operator $C^T$ defined
by (\ref{CT_dyn_repr_JM}) is a finite-rank operator,
$C^T\mathcal{F}^T:=\mathcal{F}^T_1$, and $C^T|_{\mathcal{F}^T_1}$
is an isomorphism.
\end{theorem}
\begin{proof}
The necessity of conditions has already been shown, we need to
prove the sufficiency part. The response function having form
(\ref{resp_func_JM}) gives rise to connecting operator defined by
formula (\ref{CT_dyn_repr_JM}), or equivalently
(\ref{CT_sp_repr_JM}). The space of controls denoted again by
$\mathcal{F}^T=L_2(0,T)$, is equipped with the scalar product
$f,g\in \mathcal{F}^T$,
$\left(f,g\right)_{\mathcal{F}^T}=\int_0^Tf(t)g(t)\,dt$. We
introduce the finite-dimensional space
$\mathcal{F}^T_1:=C^T\mathcal{F}^T$, where
$\operatorname{dim}\mathcal{F}^T_1=N$. Using the invertibility of
$C^T$ in $\mathcal{F}^T_1$, we find the unique solution to the
equation
\begin{equation*}
\left(C^Tf\right)(t)=r(T-t),
\end{equation*}
and denote it by $f^1$. In accordance with Proposition
\ref{Prop3}, one has the equality
$\left(C^Tf^1,f^1\right)_{\mathcal{F}^T}=1$. We set
\begin{eqnarray*}
b_1:=-\left(\left(C^Tf^1\right)'',f^1\right)_{\mathcal{F}^T},\\
\mathcal{F}^T_1\ni h^2:=-\left(C^Tf^1\right)''-b_1C^Tf^1.
\end{eqnarray*}
Note that the definition of $h^2$ implies that
$h^2=-r''(T-t)-b_1r(T-t)\not=0.$ Then we can determine $a_1>0$ and
$f^2$ from the equalities
\begin{eqnarray*}
a_1C^Tf^2=h^2,\\
\left(C^Tf^2,f^2\right)_{\mathcal{F}^T}=1.
\end{eqnarray*}
Then we set
\begin{equation*}
b_2:=-\left(\left(C^Tf^2\right)'',f^2\right)_{\mathcal{F}^T}.
\end{equation*}
Note that the following equation holds
\begin{equation*}
-\left(C^Tf^1\right)''=b_1C^Tf^1+a_1C^Tf^{2}.
\end{equation*}
Multiplying it by $f^1$ in $\mathcal{F}^T$, we get
\begin{equation*}
a_1\left(C^Tf^{2},f^1\right)_{\mathcal{F}^T}=-b_1\left(C^Tf^1,f^1\right)_{\mathcal{F}^T}+\left(\left(C^Tf^1\right)'',f^1\right)_{\mathcal{F}^T}=0.
\end{equation*}
Introducing space $\mathcal{F}^T_{1,\,C^T}$ being
$\mathcal{F}^T_{1}$ with the norm generated by the quadratic form
$\left(C^T\cdot,\cdot\right)$, we see that $f_1\bot f_2$ in
$\mathcal{F}^T_{1,\,C^T}$. Then we set
\begin{equation*}
h^3:=-\left(C^Tf^2\right)''-a_1C^Tf^1-b_2C^Tf^2.
\end{equation*}
Using the definition of $f_2,h_2$ we can rewrite
\begin{eqnarray*}
h^3=-\frac{1}{a_1}h_2''-a_1r-\frac{b_2}{a_1}h_2
=-\frac{1}{a_1}\left(-r''(T-t)-b_1r(T-t)\right)''-a_1r(T-t)\\
-\frac{b_2}{a_1}\left(-r''(T-t)-b_1r(T-t)\right)\not=0.
\end{eqnarray*}
The function $f^3\in \mathcal{F}^T_1$ and coefficient $a_2$ are
defined from the relations
\begin{equation*}
\begin{cases}
a_2C^Tf^3=h^3,\\
\left(C^Tf^3,f^3\right)=1.
\end{cases}
\end{equation*}
Multiplying the equation
\begin{equation}
\label{eq2} -\left(C^Tf^2\right)''=a_1C^Tf^1+b_2C^Tf^2+a_2C^Tf^{3}
\end{equation}
by $f^2$ in  $\mathcal{F}^T$, we get the equality
\begin{equation*}
a_2\left(C^Tf^{3},f^2\right)_{\mathcal{F}^T}=-\left(\left(C^Tf^2\right)'',f^2\right)_{\mathcal{F}^T}-a_1\left(C^Tf^1,f^2\right)_{\mathcal{F}^T}-b_2\left(C^Tf^2,f^2\right)_{\mathcal{F}^T}=0.
\end{equation*}
Multiplying (\ref{eq2}) by $f^1$ in $\mathcal{F}^T$ we obtain the
equality
\begin{eqnarray*}
a_2\left(C^Tf^{3},f^1\right)_{\mathcal{F}^T}=-\left(\left(C^Tf^2\right)'',f^1\right)_{\mathcal{F}^T}-a_1\left(C^Tf^1,f^1\right)_{\mathcal{F}^T}-b_2\left(C^Tf^2,f^1\right)_{\mathcal{F}^T}\\
=\left(\left(C^Tf^1\right)'',f^2\right)_{\mathcal{F}^T}-a_1=0,
\end{eqnarray*}
where we used Proposition \ref{Prop2}. Thus, $f^3\bot f^2,$
$f^3\bot f^1$ in $\mathcal{F}^T_{1,\,C^T}$. We continue this
procedure and construct $f_1,\ldots,f_{N-1}$,
$a_1,\ldots,a_{N-2}$, $b_1,\ldots,b_{N-1}$, then we set
\begin{equation*}
h^N:=-\left(C^Tf^{N-1}\right)''-a_{N-2}C^Tf^{N-2}-b_{N-1}C^Tf^{N-1},
\end{equation*}
We find $f_N$ and $a_{N-1}$ from the conditions
\begin{equation*}
\begin{cases}
a_{N-1}C^Tf^N=h^N,\\
\left(C^Tf^N,f^N\right)_{\mathcal{F}^T}=1,
\end{cases}
\end{equation*}
set
$b_N:=-\left(\left(C^Tf^N\right)'',f^N\right)_{\mathcal{F}^T}$,
and as it was described, show that $f^N$ is orthogonal to
$f^1,\ldots,f^{N-1}$ in $\mathcal{F}^T_{1,\,C^T}$. Then we define
\begin{equation*}
h^{N+1}:=-\left(C^Tf^{N}\right)''-a_{N-1}C^Tf^{N-1}-b_{N}C^Tf^{N},
\end{equation*}
Repeating the described procedure, we look for
\begin{equation}
\begin{cases}
a_{N}C^Tf^{N+1}=h^{N+1},\\
\left(C^Tf^{N+1},f^{N+1}\right)_{\mathcal{F}^T}=1
\end{cases}
\end{equation}
Then the function $f^{N+1}$ constricted in such a way, would be
orthogonal to $f^1,\ldots,f^N$ in $\mathcal{F}^T_{1,\,C^T}$. But
since $\operatorname{dim}\mathcal{F}^T_{1,C^T}=N$, we have that
$f^{N+1}=0$ and thus $h^{N+1}=0.$ All aforesaid leads to the
following relations for $f^1,\ldots,f^N$:
\begin{equation}
\label{CT_syst} -\begin{pmatrix} C^Tf^1 \\ \cdot \\ C^Tf^k \\ \cdot \\
C^Tf^{N}
\end{pmatrix}''=A \begin{pmatrix}
C^Tf^1 \\ \cdot \\ C^Tf^k \\ \cdot \\ C^Tf^{N}
\end{pmatrix}, \quad
A=\begin{pmatrix} b_1 & a_1 & 0 & \cdot & 0 \\
a_1 & b_2 & a_2 & \cdot & 0\\
\cdot & \cdot& \cdot& \cdot & 0\\
0 &\cdot & a_{N-2} & b_{N-1} & a_{N-1}\\
 0 &\cdot & 0& a_{N-1} & b_N
\end{pmatrix}.
\end{equation}
For the matrix $A$ which has been recovered by above procedure, we
consider the dynamical system with control
$F(t)=(f(t),0,\ldots,0)$, $f\in L_2(0,T)$ on the same time
interval $0<t<T$:
\begin{equation}
\label{DnSst_til}
\begin{cases}
\widetilde u_{tt}(t)-A\widetilde u(t)=F(t),\quad 0<t<T,\\
\widetilde u(0)=\widetilde u_t(0)=0,
\end{cases}
\end{equation}
whose response function (\ref{Resp_def}), (\ref{resp_func_JM}) we
denote by $\widetilde r$. It has the representation
\begin{equation}
\label{resp_func_JM_til} \widetilde
r(t)=\sum_{k=1}^{N}\frac{1}{\widetilde\rho_k}\widetilde
S_k(t),\quad \sum_{k=1}^N\frac{1}{\widetilde\rho_k}=1,
\end{equation}
where
\begin{equation*}
\widetilde S_k(t)=\begin{cases}
\frac{\sin{\sqrt{\widetilde\lambda_k}t}}{\sqrt{\widetilde\lambda_k}},\quad
\widetilde\lambda_k>0,\\
\frac{\operatorname{sh}{\sqrt{\widetilde
|\lambda_k|}t}}{\sqrt{\widetilde |\lambda_k|}},\quad
\widetilde\lambda_k<0,\\
t,\quad \widetilde\lambda_k=0,
\end{cases}
\end{equation*}
Our aim is to show that $r=\widetilde r$. The connecting operator
$\widetilde C^T$ for (\ref{DnSst_til}) admits representations
(\ref{CT_dyn_repr_JM}), (\ref{CT_sp_repr_JM}) with
$r,S_k,\rho_k,\lambda_k$ substituted by ones with tilde. We remind
definitions:
\begin{eqnarray*}
A\widetilde \varphi_k=\widetilde \lambda_k
\widetilde \varphi_k,\quad \widetilde\varphi_k^1=1,\\
\widetilde\rho_k=\left(\widetilde\varphi_k,\widetilde\varphi_k\right).
\end{eqnarray*}

Note that the family of special controls $\widetilde f^1,\ldots
\widetilde f^N$ (see subsection \ref{Kr_eq}) for the system
(\ref{DnSst_til}) satisfy the system
\begin{eqnarray*}
-\begin{pmatrix} \widetilde C^T \widetilde f^1 \\ \cdot \\ \widetilde C^T \widetilde f^k \\ \cdot \\
\widetilde C^T \widetilde f^{N}
\end{pmatrix}''=A \begin{pmatrix}
\widetilde C^T\widetilde f^1 \\ \cdot \\ \widetilde C^T \widetilde f^k \\ \cdot \\
\widetilde C^T \widetilde f^{N}
\end{pmatrix},
\end{eqnarray*}
with the same matrix $A$. As functions $C^Tf^1,\ldots,C^Tf^N\in
\mathcal{F}^T_1$ and $\widetilde C^T\widetilde f^1,\ldots
\widetilde C^T\widetilde f^N\in \widetilde{\mathcal{F}}^T_1$
satisfy the same second order system, we have that
$\widetilde{\mathcal{F}}^T_1=\mathcal{F}^T_1$. Or in other words
in representations for $r$, $\widetilde r$ we have that
$\lambda_k=\widetilde \lambda_k$, $k=1,\ldots,N$, and consequently
$S_k=\widetilde S_k$.

The function $f^1$ satisfies the equation $C^Tf^1=r(T-t)$, using
the representations (\ref{CT_sp_repr_JM}) and (\ref{resp_func_JM})
we have that
\begin{equation*}
C^Tf^1=\sum_{k=1}^N\frac{1}{\rho_k}\int_0^TS_k(T-t)S_k(T-s)f^1(s)\,ds=\sum_{k=1}^N\frac{1}{\rho_k}S_k(T-t),
\end{equation*}
which yields equalities
\begin{equation*}
\int_0^TS_k(T-\tau)f^1(\tau)\,d\tau=1,\quad k=1,\ldots,N.
\end{equation*}
Since for $\widetilde f^1$ the equation $\widetilde C^T\widetilde
f^1=\widetilde r(T-t)$ holds, the same arguments lead to
equalities
\begin{equation*}
\int_0^T\widetilde S_k(T-\tau)\widetilde f^1(\tau)\,d\tau=1,\quad
k=1,\ldots,N.
\end{equation*}
But since $\widetilde S_k(T-\tau) =S_k(T-\tau)$, $k=1,\ldots,N,$
we have that $\widetilde f^1=f^1$.

In view of Proposition \ref{Prop_contr}, the components of
eigenvectors $\widetilde \varphi_n$ have the following form:
\begin{equation*}
\widetilde{\varphi}_n^k=\int_0^T S_n(T-t)\widetilde
f^k(t)\,dt,\quad k=1,\ldots,N
\end{equation*}
here we use that $\widetilde S_k(T-\tau) =S_k(T-\tau)$,
$k=1,\ldots,N$. Introducing the vectors $\psi_n$ whose components
are defined by the rule
\begin{equation*}
\psi_n^k=\int_0^T S_n(T-t) f^k(t)\,dt, \quad k=1,\ldots,N,
\end{equation*}
we will show that $\psi_n$ are also eigenvectors of $A$
corresponding to $\lambda_n$. To do so we observe that
\begin{eqnarray*}
C^Tf^k=\sum_{n=1}^N\frac{1}{\rho_n} S_n(T-t)\psi_n^k,\quad k=1,\ldots,N,\\
-\left(C^Tf^k\right)''=\sum_{n=1}^N\lambda_n\frac{S_n(T-t)}{\rho_n}\psi_n^k,\quad
k=1,\ldots,N.
\end{eqnarray*}
Plugging the above quantities into (\ref{CT_syst}) and equating
coefficients at $S_n(T-t)$, we see that $A\psi_n=\lambda_n\psi_n$.
The equality $\widetilde f^1=f^1$ implies the fact that first
components of eigenvectors $\widetilde \varphi_n$ and $\psi_n$
coincide:
\begin{equation*}
\widetilde \varphi_n^1=\int_0^TS_n(T-t)\widetilde
f^1(t)\,dt=\int_0^TS_n(T-t) f^1(t)\,dt=\psi_n^1.
\end{equation*}
But then we have that
\begin{equation*}
\widetilde \psi_n=\varphi_n,\quad \widetilde f^k=f^k, \quad
n=1,\ldots,N.
\end{equation*}
The following equations hold (we count that $\widetilde f_k=f_k$):
\begin{eqnarray*}
\left(C^Tf^1\right)(t)=r(T-t),\\
\left(\widetilde{C}^T f^1\right)(t)=\widetilde r(T-t)
\end{eqnarray*}
Evaluating scalar products we see that
\begin{equation*}
(r,f^1)_{\mathcal{F}^T}=\left(C^Tf^1,f^1\right)_{\mathcal{F}^T}=1,\quad
(r,f^k)_{\mathcal{F}^T}=\left(C^Tf^1,f^k\right)_{\mathcal{F}^T}=0,\quad
k=2,\ldots,N
\end{equation*}
by construction. And
\begin{equation*}
(\widetilde r,f^1)_{\mathcal{F}^T}=\left(\widetilde C^T\widetilde
f^1,\widetilde f^1\right)_{\mathcal{F}^T}=1,\quad (\widetilde
r,f^k)_{\mathcal{F}^T}=\left(\widetilde C^T\widetilde
f^1,\widetilde f^k\right)_{\mathcal{F}^T}=0,\quad k=2,\ldots,N
\end{equation*}
by properties of the operator $\widetilde C^T.$ All aforesaid
implies $r=\widetilde r$, which completes the proof.

\end{proof}

\section{IP for a finite Krein--Stieltjes string.}

\subsection{Forward problem.}
In this section we consider the special case of a dynamic IP for a
string whose mass is a nondecreasing bounded function, see
\cite{Kr}. The authors in \cite{MM3,MM31} set up the dynamic
forward and inverse problems for such a string, here we would like
to elaborate and improve some of results obtained in these papers.

Let on an interval $(0,l)$ a nondecreasing bounded function (a
mass of a string) $M(x)$: $M(0)<M(x)<M(l)$ be given.
Following \cite{KK,DMcK}, we introduce the domain consisting of
\emph{continued functions}
\begin{eqnarray*}
D_M:=\left\{\left[u(x),u'_-0),u'_+(l)\right]\,\bigl|\,
u(x)=a+bx-\int_0^x(x-s)g(s)M'(s)\,ds;\right.\\
\left.u'_-(0)=b;\,u'_+(l)=b-\int_0^lg(s)M'(s)\,ds \right\},
\end{eqnarray*}
where $u'_-(0),u'_+(l)$ are left and right derivatives, $g$ is
$M-$summable; and define the \emph{generalized differential
operation} $l_M[u]$ on $D_M$ by the rule
\begin{equation*}
\label{L_oper} l_M[u]=g(x),\quad \text{$M-$almost everywhere}.
\end{equation*}
Note that if $M$ is $C^1$ smooth, then $
l_M[y]=-\frac{1}{M'(x)}\frac{d^2y(x)}{dx^2}. $

We fix $T>0$ and consider the following initial-boundary value
problem (IBVP):
\begin{equation}
\label{L_eq}
\begin{cases}
u_{tt}(x,t)+l_M[u]=0,\quad 0<x<l,\, 0<t<T,\\
u(x,0)=0,\, u(x,l)=0,\quad 0\leqslant x\leqslant l,\\
u(0,t)=f(t),\quad 0\leqslant t\leqslant T,
\end{cases}
\end{equation}
where the \emph{boundary control} $f\in L_2(0,T)$. If $M$ is
$C^1-$smooth, $M'(x)\geqslant \delta>0$, $0\leqslant x\leqslant
l$, then (\ref{L_eq}) corresponds to the IBVP for a wave equation
\cite{BM_1,BL}. We assume that the string is a Stieltjes-Krein
string, i.e. the mass $M$ is a piecewise constant function: let
$0=x_0<x_1<x_2<\ldots<x_{N-1}<x_N<x_{N+1}=l$, $m_i>0,$
$i=1,\ldots,N$, $l_i=x_i-x_{i-1}$, $i=1,\ldots N+1,$ and the
density $dM$ has a form $dM(x)=\sum_{i=1}^{N}m_i\delta(x-x_i)$. In
this specific case we have that
\begin{equation*}
l_M[u](x)=\frac{u'(x+0)-u'(x-0)}{m_x},\quad m_x:=\begin{cases} M(+0)-M(0),\quad x=0,\\
M(x+0)-M(x-0),\quad 0<x<l,\\
M(l)-M(l-0),\quad x=l.
\end{cases}
\end{equation*}
and thus the IBVP (\ref{L_eq}) is equivalent to the IBVP for the
vector function (we keep the same notation)
$u(t)=\left(u_1(t),\ldots,u_{N-1}(t)\right)$:
\begin{equation}
\label{DynSyst}
\begin{cases}
Mu_{tt}=Au+\widetilde{f},\quad t> 0,\\
u(0)=0,\, u_t(0)=0,
\end{cases}
\end{equation}
where
\begin{eqnarray*}
A=\begin{pmatrix} b_1 & a_1&0&   \ldots & 0\\
a_1& b_2 & a_2& \ldots & 0\\
\ldots & \ldots& \ldots & \ldots \\
0 & 0& \ldots&  b_{N-1} &a_{N-1}\\
0 & 0& \ldots& a_{N-1} &b_{N}
\end{pmatrix},\,
M=\begin{pmatrix} m_1 & 0&  \ldots & 0\\
0& m_2 &  \ldots & 0\\
\ldots & \ldots& \ldots & \ldots \\
0 & \ldots& 0 & m_{N}
\end{pmatrix}, \,
\widetilde f=\begin{pmatrix} \frac{f}{l_1}\\0\\
\ldots\\0\end{pmatrix},
\end{eqnarray*}
and the entries of the matrix $A$ are given by
\begin{equation}
\label{coeff} a_i=\frac{1}{l_{i+1}},\,\,
b_i=-\frac{l_i+l_{i+1}}{l_il_{i+1}}.
\end{equation}
The solution to (\ref{DynSyst}) is denoted by $u^f$. The dynamic
IP consists in recovering the string, i.e. the coefficients
$l_N,m_i,l_i$, $i=1,\ldots,N$ or, what is equivalent, the entries
of the matrices $A$ and $M$, from the knowledge of a
\emph{response operator} $R^T,$ defined by the rule
\begin{equation}
\label{Resp}
\left(R^Tf\right)(t):=u^f_1(t),\quad 0<t<T.
\end{equation}


Let us denote $A_M=M^{-\frac{1}{2}}AM^{-\frac{1}{2}}$,
and let $C$ be unitary matrix such that
$C^*A_MC=\operatorname{diag}\{\lambda_1,\ldots,\lambda_{N}\}$,
$\lambda_1<\lambda_2<\ldots<\lambda_{N}$; let
$A_Mg_k=\lambda_kg_k,$ $\|g_k\|=1$. Note that all eigenvalues of
$A_M$ are negative: it is a characteristic properties of a
Stieltjes-Krein string, see \cite{EK}. 
\begin{lemma}
The solution to (\ref{DynSyst}) admits the representation
\begin{equation}
\label{Sol_dyn_repr} u_j^f(t)=\frac{1}{l_1}\sum_{k=1}^{N}f_k(t)
c_{k1}c_{kj},\quad j=1,\ldots,N,
\end{equation}
where
\begin{equation*}
c_{kj}=\frac{g_{kj}}{\sqrt{m_j}},\quad f_k(t)=\int_0^t
f(\tau)\frac{\sin{\sqrt{|\lambda_k|}(t-\tau)}}{\sqrt{|\lambda_k|}}\,d\tau,\quad
k,j=1,\ldots,N.
\end{equation*}
\end{lemma}
The following Cauchy problem for the difference equation
\begin{equation}\label{CPfDE}
\begin{cases}
a_n\phi_{n+1}+a_{n-1}\phi_{n-1}+b_n\phi_n=\lambda m_n\phi_n,\quad
n=1,\ldots,N,\\
\phi_0=0,\phi_1=1,
\end{cases}
\end{equation}
determines the set of polynomials
$\{0,1,\phi_2(\lambda),\ldots,\phi_N(\lambda),\phi_{N+1}(\lambda)\}$.
Note that $\lambda_1,\ldots,\lambda_{N}$ are roots of the equation
$\phi_{N+1}(\lambda)=0$. Upon introducing the vectors and the
coefficients
\begin{equation*}
\varphi(\lambda)=\begin{pmatrix} \phi_1(\lambda)\\
\cdot\\ \cdot\\ \phi_{N}(\lambda)\end{pmatrix},\quad\varphi_k=\begin{pmatrix} \phi_1(\lambda_k)\\
\cdot\\ \cdot\\ \phi_{N}(\lambda_k)\end{pmatrix},\quad
\rho_k=\left(M\varphi_k,\varphi_k\right).
\end{equation*}
we call by \emph{spectral data} and spectral function $\rho$ the
following objects:
\begin{equation}
\label{measure} \left\{\lambda_i,\rho_i\right\}_{i=1}^{N},\quad
\rho(\lambda)=\sum_{\{k:\lambda_k<\lambda\}}\frac{1}{\rho_k}.
\end{equation}
The standard application of a Fourier method (see
\cite{AI,AM,MM1}) yields
\begin{lemma}
The solution to (\ref{DynSyst}) admits the spectral representation
\begin{eqnarray}
\label{Sol_spec_repr} u^f(t)=\sum_{k=1}^{N}h_k(t)\varphi_k, \quad
h_k(t)=\frac{1}{l_1\rho_k}\int_0^t
f(\tau)\frac{\sin{\sqrt{|\lambda_k|}(t-\tau)}}{\sqrt{|\lambda_k|}}\,d\tau,\\
u^f(t)=\frac{1}{l_1}\int_{-\infty}^\infty\int_0^t\frac{\sin{\sqrt{|\lambda_k|}(t-\tau)}}{\sqrt{|\lambda_k|}}f(\tau)\,d\tau\varphi(\lambda)\,d\rho(\lambda).\notag
\end{eqnarray}
\end{lemma}

\subsection{Operators of the Boundary control method.}

\subsubsection{Response operator. } The \emph{outer space} of the
system (\ref{DynSyst}), the space of controls is denoted by
$\mathcal{F}^T:=L_2(0,T)$. The \emph{response operator} $R^T:
\mathcal{F}^T\mapsto \mathcal{F}^T$ is introduced by (\ref{Resp}).
Making use of (\ref{Sol_spec_repr}) implies the representation for
$R^T$:
\begin{eqnarray}
\left(R^Tf\right)(t)=u_1^f(t)=\sum_{k=1}^{N}h_k(t)=\int_0^t
r(t-s)f(s)\,ds,\notag\\
\text{where}\quad
r(t)=\frac{1}{l_1}\sum_{k=1}^{N}\frac{\sin{\sqrt{|\lambda_k|}t}}{\sqrt{|\lambda_k|}\rho_k}\,\,\text{is
a \emph{response function}}. \label{resp_func}
\end{eqnarray}

The \emph{inner space}, i.e. the space of states of
(\ref{DynSyst}) is denoted by $\mathcal{H}^N:=\mathbb{R}^{N}$, so
for $T>0$, $u^f(T)\in \mathcal{H}^T$. The metric in
$\mathcal{H}^N$ is given by
\begin{equation*}
(a,b)_{\mathcal{H}^T}=\left(Ma,b\right)_{\mathbb{R}^{N}},\quad
a,b\in\mathcal{H}^N.
\end{equation*}
The \emph{control operator} $W^T: \mathcal{F}^T\mapsto
\mathcal{H}^N$ is introduced by the rule:
\begin{equation*}
W^Tf=u^f(T).
\end{equation*}
Due to (\ref{Sol_spec_repr}) we have that
$W^Tf=\sum_{k=1}^{N}h_k(T)\varphi_k$. 
For a prescribed state $a\in
\mathcal{H}^N$, $a=\sum_{k=1}^{N-1}a_k\varphi_k$ the existence of
a control $f\in \mathcal{F}^T$ such that $W^Tf=a$ is equivalent
(see (\ref{Sol_spec_repr})) to the trigonometric moment problem:
to find $f\in \mathcal{F}^T$ such that
\begin{equation*}
a_k=h_k(T)=\frac{1}{l_1\rho_k}\int_0^T
f(\tau)\frac{\sin{\sqrt{|\lambda_k|}(T-\tau)}}{\sqrt{|\lambda_k|}}\,d\tau,\quad
k=1,\ldots,N.
\end{equation*}
Clearly, the solution to this moment problem is not unique. We
introduce a space
\begin{equation*}
\mathcal{F}^T_1={\operatorname{span}\left\{\frac{\sin{\sqrt{|\lambda_k|}(T-t)}}{\sqrt{|\lambda_k|}}\right\}}_{k=1}^N.
\end{equation*}
The following lemma answers the question on the controllability of
(\ref{DynSyst}):
\begin{lemma}
\label{LemmaCont} The operator ${W}^T$ is an isomorphism between
$\mathcal{F}^T_1$ and $\mathcal{H}^N$.
\end{lemma}


The \emph{connecting operator} $C^T:\mathcal{F}^T\mapsto
\mathcal{F}^T$ is defined by $C^T:=\left(W^T\right)^*W^T$, so by
the definition for $f,g\in \mathcal{F}^T$ one has
\begin{equation}
\label{C_T}
\left(C^Tf,g\right)_{\mathcal{F}^T}=\left(u^f(T),u^g(T)\right)_{\mathcal{H}^N}=\left(W^Tf,W^Tg\right)_{\mathcal{H}^N}.
\end{equation}
It is crucial in the BC method that $C^T$ can be expressed in
terms of inverse data:
\begin{theorem}
The connecting operator admits the representation in terms of
dynamic inverse data:
\begin{equation}
\label{CT_dyn_repr}
\left(C^Tf\right)(t)=\frac{1}{l_1}\int_0^T\int_{|t-s|}^{2T-s-t}r(\tau)\,d\tau
f(s)\,ds,
\end{equation}
and in terms of spectral inverse data:
\begin{equation}
\label{CT_sp_repr} \left(C^Tf\right)(t)=\int_0^T
\frac{1}{l_1}\sum_{k=1}^{N}\frac{\sin{\sqrt{|\lambda_k|}(T-t)}\sin{\sqrt{|\lambda_k|}(T-s)}}{{|\lambda_k|}\rho_k}f(s)\,ds,
\end{equation}
\end{theorem}
\begin{proof}
Introducing the function $\psi(t,s):=\left(Mu^f(t),u^g(s)\right)$,
where $f,g\in \mathcal{F}^T$, one can show that $\psi$ satisfies
the equation
\begin{equation*}
\psi_{tt}-\psi_{ss}=\frac{1}{l_1}\left(f(t)(Rg)(t)-g(s)(Rf)(t)\right).
\end{equation*}
Solving it by the d'Alembert method and noting that
$\psi(T,T):=\left(C^Tf,g\right)_{\mathcal{F}^T}$, yields
(\ref{CT_dyn_repr}). The formula (\ref{CT_sp_repr}) follows from
(\ref{C_T}) and (\ref{Sol_spec_repr}).
\end{proof}

\begin{remark}
The formula (\ref{CT_sp_repr}) implies that
$\mathcal{F}^T_1=C^T\mathcal{F}^T$, that is, $\mathcal{F}^T_1$ is
determined by inverse data. The important properties of functions
from $\mathcal{F}^T_1$ is that
\begin{equation}
\label{FT_proper} f(T)=0,\quad \text{for all}\,\, f\in
\mathcal{F}^T_1.
\end{equation}
\end{remark}

\subsection{Inverse problem.}

In  this section we outline the method of solving the IP for the
system (\ref{DynSyst}) based on Krein equations.

By $f_k^T\in \mathcal{F}^T_1$ we denote controls, driving the
system (\ref{DynSyst}) to prescribed \emph{special states}
\begin{equation*}
d_k\in \mathcal{H}^N,
\,d_k=\left(0,\ldots,\frac{1}{m_k},\ldots,0\right),\quad
k=1,\ldots,N.
\end{equation*}
In other words, $f^T_k$ are solutions to the equations
$W^Tf_k^T=d_k$, $k=1\,\ldots,N$. By Lemma \ref{LemmaCont} we know
that such controls exist and are unique. It is important that they
can be found as a solution to a Krein equation.

\begin{theorem} The control $f_1^T$ can be found as a solution to
the following equation:
\begin{equation}
\label{Cont_f1} \left(C^Tf_1^T\right)(t)=r(T-t),\quad 0<t<T.
\end{equation}
The controls $f^T_k$ satisfy the system:
\begin{equation}
\label{Cont_system}
\begin{cases}
m_1\left(C^Tf_1^T\right)''=b_1C^Tf^T_1+a_1C^Tf^T_{2},\\
m_k\left(C^Tf_k^T\right)''=a_{k-1}C^Tf^T_{k-1}+b_kC^Tf^T_k+a_kC^Tf^T_{k+1},\quad
k=2,\ldots,N-1,\\
m_N\left(C^Tf_N^T\right)''=a_{N-1}C^Tf^T_{N-1}+b_NC^Tf^T_N,
\end{cases}
\end{equation}
\end{theorem}
\begin{proof}
Let us check (\ref{Cont_f1}): for $g\in \mathcal{F}^T$ we have:
\begin{equation*}
\left(C^Tf_1^T,g\right)_{\mathcal{F}^T}=\left(Mu^{f_1^T}(T),u^g(T)\right)_{\mathcal{H}^N}=u_1^g(T)=\int_0^Tr(T-\tau)g(\tau)\,d\tau,
\end{equation*}
which, due to the arbitrariness of $g$, yields (\ref{Cont_f1}).

From (\ref{Sol_spec_repr}) we see that
\begin{equation*}
 u^{f_k^T}_n(T)=\frac{1}{l_1}\sum_{i=1}^{N-1} \frac{1}{\rho_i}\int_0^T
 f^T_k(\tau)\frac{\sin{\sqrt{|\lambda_i|}(T-\tau)}}{\sqrt{|\lambda_i|}}\,d\tau \varphi^n_i=\frac{\delta_{kn}}{m_k}.
 \end{equation*}
Multiplying the above equality by $\varphi^n_j m_n$ and taking the
sum with respect to $n$ we get that
\begin{equation}\label{coef_eq}
\frac{1}{l_1} \int_0^T
f^T_k(\tau)\frac{\sin{\sqrt{|\lambda_j|}(T-\tau)}}{\sqrt{|\lambda_j|}}\,d\tau = \varphi^k_j.
\end{equation}
Now we can rewrite (\ref{CT_sp_repr}) for $f=f^T_k$:
\begin{equation}\label{Ctf}
(C^Tf_k^T)(t)= \sum_{i=1}^{N-1}
\varphi^k_i\frac{\sin{\sqrt{|\lambda_i|}(T-t)}}{\sqrt{|\lambda_i|}\rho_i}.
\end{equation}
If we take the second derivative of the last equality and multiply
it by $m_k$, we get that
\begin{equation}
\label{CK_sum} m_k(C^Tf_k^T)''(t)= \sum_{i=1}^{N-1} m_k
\lambda_i\varphi^k_i\frac{\sin{\sqrt{|\lambda_i|}(T-t)}}{\sqrt{|\lambda_i|}\rho_i}.
\end{equation}
We remind that $\varphi_i$ are eigenvectors of the problem
$A\varphi_i=\lambda_iM\varphi_i$ (see (\ref{CPfDE})), so we plug
$$
m_k \lambda_i\varphi^k_i = a_k\varphi^{k+1}_i+a_{k-1}\varphi^{k-1}_i+b_k\varphi^k_i
$$
into (\ref{CK_sum}) and obtain the following
\begin{multline*}
m_k(C^Tf_k^T)''(t)= \sum_{i=1}^{N-1}  (a_k\varphi^{k+1}_i+a_{k-1}\varphi^{k-1}_i+b_k\varphi^k_i)\frac{\sin{\sqrt{|\lambda_i|}(T-t)}}{\sqrt{|\lambda_i|}\rho_i}=\\
=a_{k-1}C^Tf^T_{k-1}+b_kC^Tf^T_k+a_kC^Tf^T_{k+1},
\end{multline*}
where the last equality is due to (\ref{Ctf}). This completes the
proof of the statement.
\end{proof}

\begin{remark}
Parameters of the string can be recovered by
\begin{eqnarray}
l_1=-\frac{\|f_1^T\|^2_{L_2(0,T)}}{\left(f_1^T\right)'(T)},\quad
m_k=\frac{1}{\left(C^Tf_k^T,f_k^T\right)_{\mathcal{F}^T}},\quad
k=1,\ldots,N,\label{mk}\\
b_k=m_k^2\left(\left(C^Tf_k^T\right)'',f_k^T\right)_{\mathcal{F}^T},\quad
k=1,\ldots,N.\label{bk}
\end{eqnarray}
\end{remark}
Indeed, the formula for $m_k$ is a consequence of a definition of
controls $f_k^T$, the formula for $b_k$ is a consequence of
equation (\ref{Cont_system}) and orthogonality of controls
$f_k^T$.  Formula for $l_1$ can be obtained in the following way:
(we remind that $\varphi^1_j=1$ for $j=1,..,N$) from
(\ref{coef_eq}) we have that
 $$
 \int_0^T
 f^T_1(\tau)\frac{\sin{\sqrt{|\lambda_j|}(T-\tau)}}{\sqrt{|\lambda_j|}}\,d\tau = l_1.
 $$
We know that  $f_1^T\in\mathcal{F}^T_1$. It means that for some
$d_k\in \mathbb{R}$ the control has a form
$$
f^T_1(\tau)=\sum_{k=1}^N d_k\frac{\sin{\sqrt{|\lambda_j|}(T-\tau)}}{\sqrt{|\lambda_j|}}.
$$
Computing the $L_2$-norm of the control we obtain that
$$
\|f_1^T\|^2_{L_2}= (f_1^T(\cdot),f_1^T(\cdot))=\int_{0}^{T}
f_1^T(\tau)\sum_{k=1}^N
d_k\frac{\sin{\sqrt{|\lambda_j|}(T-\tau)}}{\sqrt{|\lambda_j|}}=l_1\sum_{k=1}^{N}
d_k,
$$
and computing the first derivative of $f_1^T(t)$ at the point
$t=T$ gives
$$
(f_1^T)'|_{t=T}=\sum_{k=1}^{N} d_k
\cos{\sqrt{|\lambda_j|}(T-t)}(-1)|_{t=T}= -\sum_{k=1}^{N} d_k
$$
From the last two equations we deduce the formula for $l_1$.


\begin{theorem}
The function $r$ is a response function for the dynamical system
(\ref{L_eq}) if and only if it has a form (\ref{resp_func}),
operator $C^T$ defined by (\ref{CT_dyn_repr}) is a finite-rank
operator, $C^T\mathcal{F}^T:=\mathcal{F}^T_1$, and
$C^T|_{\mathcal{F}^T_1}$ is an isomorphism.
\end{theorem}
\begin{proof}
The necessity of conditions has already been shown. We need to
prove the sufficiency part. 
The proof repeats arguments of Theorem \ref{thm3}. If we have
function $r$ of the form (\ref{resp_func}) and we know that $C^T$
constructed by is an isomorphism on its image, we can uniquely
solve equation (\ref{Cont_f1}) to construct $f_1$. Then we define
$l_1$, $m_1$, $b_1$ by (\ref{mk}), (\ref{bk}). After that we
construct $a_1$ from (\ref{coeff}). To construct $f_2$ we consider
(\ref{Cont_system}) and see that all terms are defined except
$C^Tf_2$. Therefore $C^Tf_2$ and $f_2$ are known. After that we
define  $m_2$, $b_2$, $a_2$ and etc. Using (\ref{Cont_system}) we
see that $(C^Tf_k,f_l)=\delta_{kl}/m_k$.

Then, using coefficients $a_k$, $b_k$, $m_k$, $k=1,\ldots,N$ we
construct system of the form (\ref{L_eq}). It remains to show that
the response function of this system is exactly $r$. For a moment
we call this function by $\widetilde r$ and connecting operator
for this system by $\widetilde C^T$. Then $\widetilde r$ has the
form
\begin{equation*}
\widetilde r(t)=\frac{1}{\widetilde l_1}
\sum_{k=1}^{N-1}\frac{\sin{\sqrt{|\widetilde\lambda_k|}t}}{\sqrt{|\widetilde\lambda_k|}.
\widetilde\rho_k}
\end{equation*}
Note, that $\widetilde\lambda_k$ is an eigenvalue of matrix $A_M =
M^{-\frac{1}{2}} A M^{-\frac{1}{2}} $. At the same time from
(\ref{Cont_system}) we see that $\lambda_k$ is an eigenvalue of
the same matrix. Therefore $\widetilde\lambda_k=\lambda_k$. And
from (\ref{coef_eq}) we deduce that $\tilde f_k =f_k$ and
therefore $\rho_k=\widetilde\rho_k$ and $\widetilde r(t)= r(t)$.
\end{proof}

\noindent{\bf Acknowledgments}

The research of Victor Mikhaylov was supported in part by RFBR
18-01-00269, 20-01-00627 and by the Ministry of Education and
Science of Republic of Kazakhstan under grant AP05136197. Alexandr
Mikhaylov was supported by RFBR 18-01-00269.  Alexandr Mikhaylov
and Victor Mikhaylov were supported by Volkswagen Foundation
project "From Modeling and Analysis to Approximation". The authors
would like to express gratitude to anonymous referee for valuable
comments and remarks.

\end{document}